\documentclass[12pt, leqno, oneside]{amsart}
\usepackage{geometry}
\usepackage{enumerate}
\newtheorem{theorem}{Theorem}
\newtheorem{corollary}[theorem]{Corollary}
\newtheorem{lemma}[theorem]{Lemma}
\newtheorem{proposition}[theorem]{Proposition}
\theoremstyle{definition}

\theoremstyle{remark}
\newtheorem{remark}[theorem]{Remark}
\def\prop{\operatorname{Prop}}
\def\aut{\operatorname{Aut}}

\title[Proper maps]
{Proper holomorphic mapppings between Reinhardt domains in $\mathbb C^2$}
\author[\L.~Kosi\'nski]
{\L ukasz Kosi\'nski}
\address{Instytut Matematyki\\ Uniwersytet Jagiello\'nski\\ Reymonta 4
30-059 Krak\'ow\\ Poland} \email{lukasz.kosinski@im.uj.edu.pl} 

\begin{document}
\baselineskip=16pt
\begin{abstract}
\baselineskip=16pt We describe all possibilities of existence of non-elementary proper holomorphic maps between non-hyperbolic Reinhardt domains
in $\mathbb C^2$ and the corresponding pairs of domains.
\end{abstract}

\maketitle

\section{Introduction}
Given $\alpha\in\mathbb R^n,$ and $z\in\mathbb C_*^n$ we put $|z^{\alpha}|:=|z|^{\alpha_1}\ldots |z|^{\alpha_n}$ whenever it makes sense. Let
$\mathbb A_{r^-,r^+}=\{z\in\mathbb C:\ r^-<|z|<r^+\}$ for $-\infty<r^-<r^+<\infty,\ r^+>0.$ By $\mathbb D$ we always denote the unit disc in
$\mathbb C.$ For a domain $D\subset\mathbb C^n,$ $D\setminus\{0\}$ is denoted by $D_*.$

Following \cite{dis}, for $A=(A_k^j)_{j=1,\ldots,m,\ k=1,\ldots,n}\in\mathbb Z^{m\times n}$ and $b=(b_1,\ldots,b_m)\in\mathbb C_*^m$ we define:
\begin{align*}&\varphi_A(z):=z^A:=(z^{A^1},\ldots,z^{A^m}),\quad
z\in\mathbb C_*^n,\\
&\varphi_{A,b}(z):=(b_1z^{A^1},\ldots,b_mz^{A^m}),\quad z\in\mathbb C_*^n.\end{align*} Such maps are called \textit{elementary algebraic} (or
briefly \textit{elementary} maps).

The aim of this paper is to describe non-elementary proper holomorphic maps between non-hyperbolic Reinhardt domains in $\mathbb C^2$ as well as
the corresponding pairs of domains. Additionally, we  obtain some partial results for proper maps between domains of the form $\mathbb C_*^2,$
$\mathbb C^2$ and $\mathbb C\times\mathbb C_*$ and we give some more general results related to proper holomorphic mappings.

Recall that if $D,G$ are Reinhardt domains and $f:D\to G$ is a biholomorphic mapping, then $f$ can be represented as composition of automorphism
of $D$ and $G$ and an elementary mapping between these domains (see \cite{kr} and \cite{sh}). Thus, the description of non-elementary
biholomorphic mappings between Reinhardt domains reduces to the investigation of their group of automorphisms. It is a general problem of
complex geometry of Reinhardt domains considered in many papers. In \cite{shimizu} the author using group-theoretic methods investigated the
holomorphic equivalence of bounded Reinhardt domains in $\mathbb C^n$ not containing the origin and determined automorphisms of a certain class
of Reinhardt domains. Similar results were obtained by Barrett in \cite{barrett}, however his approach was analytic. The groups of automorphisms
of all bounded Reinhardt domains containing the origin were determined in \cite{sunada}. This work has been extended in \cite{kr} by dropping
the assumption that the origin belongs to domain. The situation when domains $D$ and $G$ may be unbounded were considered for example in
\cite{shi1} and \cite{prop}.

Obviously, the problem of description of proper holomorphic mappings is harder to deal with. Proper maps between non-hyperbolic, pseudoconvex
Reinhardt domains have been considered in \cite{prop} and \cite{ja}. In the bounded case partial results were obtained in \cite{bep}, \cite{ls}
and \cite{d-sp}. The final result for bounded domains in $\mathbb C^2,$ which may be viewed as the completion of this research, has been lastly
obtained in the paper \cite{isaev} of A.V.~Isaev and N.G.~Kruzhilin. The authors explicitly described all possibilities of existence of
non-elementary proper holomorphic mappings between bounded Reinhardt domains in $\mathbb C^2.$

Our results finish the problem of characterization of proper holomorphic mappings between Reinhardt domains in $\mathbb C^2.$

\section{Preliminaries and statement of results}
It is well known that for any pseudoconvex Reinhardt domain $D$ in $\mathbb C^n$ its logarithmic image $\log D$ is convex. Moreover, any proper
holomorphic mapping between domains $D_1,D_2$ in $\mathbb C^n$ can be extended to a proper map between the envelopes of holomorphy $\widehat
D_1, \widehat D_2$ of $D_1, D_2$ respectively (see~e.g.~\cite{ker}).

Let us introduce some notation. First we define \begin{align}&V_{\iota}:=\mathbb C^{\iota-1}\times\{0\} \times\mathbb C^{n-\iota}\subset\mathbb
C^n,\quad \iota=1,\ldots,n,\quad\\\nonumber &\text{and}\quad M:=\bigcup_{\iota=1}^{n} V_{\iota}.\end{align} With a given Reinhardt domain $D$ we
associate the following constants:
\begin{enumerate}
\item[$d(D)$]:=the maximal possible dimension of the linear
subspace contained in the logarithmic image of the envelope of holomorphy of $D;$
\item[$t(D)$]:=the number of $j$
such that $\widehat D\cap V_j\neq\emptyset.$
\end{enumerate}

Moreover, in the case $D\subset \mathbb C^2$ we put
\begin{enumerate}
\item[$s(D)$]:=the number of $j=1,2$ such that $V_j \cap \widehat D$ is equal
to $\mathbb C;$
\item[$s_*(D)$]:=the number of $j=1,2$ such that $V_j \cap \widehat D$ is equal
to $\mathbb C_*;$
\end{enumerate}

It turns out that the objects introduced above are invariant under proper holomorphic mappings $f:D\to G,$ where $D,G$ are Reinhardt domains in
$\mathbb C^2,$ except for the case when $\alpha\mathbb R+\beta\subset \log D$ for some $\alpha\in\mathbb Q^2,\ \beta\in\mathbb R^2.$ In
particular, we shall obtain the following

\begin{theorem}\label{niezmienniki} Let $D,G$ be Reinhardt domains in $\mathbb C^2$
such that the set of proper holomorphic mappings from $D$ onto $G$ is non-empty. Then \begin{equation}d(D)=d(G).\end{equation}

If moreover $d(D)=d(G)=0,$ then
\begin{equation}(s(D),s_*(D),t(D))=(s(G),s_*(G),t(G)).\end{equation}
\end{theorem}

Recall here that pseudoconvex Reinhardt domains which are algebraically biholomorphic to bounded Reinhardt domains have been described by
W.~Zwonek in \cite{zwarx}. This result is of key importance for our considerations so we quote it below:
\begin{theorem}\label{wynikzwonka}
Assume that $D$ is a pseudoconvex Reinhardt domain in $\mathbb C^n.$ Then the following conditions are equivalent:
\begin{enumerate}[(i)]
\item $D$ is Brody
hyperbolic, i.e. any holomorphic mapping from $\mathbb C$ to $D$ is constant,
\item\label{hyper} \begin{enumerate}[(a)]
\item $\log D$ contains no affine lines and
\item $D\cap V_j$ is
either empty or c-hyperbolic, $j=1,\ldots,n,$ (viewed as domain in $\mathbb C^{n-1}$);
\end{enumerate}
\item $D$ is algebraically biholomorphic to a bounded Reinhardt
domain, i.e. there is $A\in\mathbb Z^{n\times n},\ |\det A|=1,$ such that $\varphi_A(D)$ is bounded and $(\varphi_A)_{|D}$ is a biholomorphism
onto the image.
\end{enumerate}
\end{theorem}
Note that a Reinhardt domain $D$ in $\mathbb C^2$ satisfies the condition~(\ref{hyper}) of Theorem~\ref{wynikzwonka} if and only if
$s(D)=s_*(D)=d(D)=0.$

Let $D_1,D_2$ be Reinhardt domains in $\mathbb C^2$ and let $f:D_1\to D_2$ be a proper holomorphic mapping. Assume that $f$ is non-elementary.
Our aim is to get the explicit formulas for the mapping $f$ as well for the domains $D_1,D_2.$

In view of Theorem~\ref{wynikzwonka} we see that case $d(D_i)=s(D_i)=s_*(D_i)=0,\ i=1,2,$ has been described in \cite{isaev}. Moreover, in
\cite{prop} and \cite{ja}  the authors gave the explicit formulas for all proper holomorphic mappings $f:D_1\to D_2$ between pseudoconvex
Reinhardt domains $D_1$ and $D_2$ such that $d(D_1)=d(D_2)=1,$ that is domains whose logarithmic image is equal to a strip or a half plane. One
may apply direct and tedious calculations which allow to determine all possibilities of the form of Reinhardt domains $D'_1, D'_2$ whose
envelopes of holomorphy are equal to $D_1$ and $D_2$ respectively and such that the restriction $f|D'_1:D'_1\to D'_2$ is proper.

On the other hand there is no proper holomorphic mapping between hyperbolic and non-hyperbolic domains (see Lemma~\ref{brod}) so we shall focus
our considerations on proper holomorphic mappings between non-hyperbolic domains.

Summing up, to obtain a desired description of the set of non-elementary proper holomorphic mappings between non-hyperbolic Reinhardt domains
$D_1,D_2$ in $\mathbb C^2$ whose envelopes of holomorphy do not contain $\mathbb C^2_*$ it suffices to confine ourselves to the cases when
$d(D_1)=d(D_2)=0$ and $s(D_1)=s(D_2)\neq0$ or $s_*(D_1)=s_*(D_2)\neq 0.$

Now we are in position to formulate the main result of this paper:
\begin{theorem}\label{main} Let $D_1,D_2$ be non-hyperbolic Reinhardt domains in $\mathbb C^2$ such
that $d(D_1)=d(D_2)=0$ and $s(D_i)\neq0$ or $s_*(D_i)\neq 0,\ i=1,2.$ Assume that there is a proper, non-elementary holomorphic mapping
$f:D_1\to D_2.$

Then the following two scenarios obtain:
\begin{enumerate}[(i)]
\item Up to a permutation of the components of $f$ and the variables,
the map $f$ has the form
\begin{equation}f(z,w)=(\mu_1 z^k B(C_1 z^{p_1}w^{q_1}),\mu_2 w^l),
\end{equation}
where $k,l\in\mathbb N,$ $p_1,q_1>0$ are relatively prime integers, $B$ is a non-constant finite Blaschke product non-vanishing at $0,$ $C_1>0$
and $\mu_1,\mu_2\in\mathbb C_*.$ In this case, the domains $D_1$ and $D_2$ have the form:
\begin{align} D_i=\{(z,w)\in\mathbb C^2:\ C_i|z|^{p_i}|w|^{q_i}<1,\
|w|<E_i\}\setminus (P_i\times\{0\}),\ i=1,2,
\end{align} where $E_1,E_2>0,\ p_2,q_2>0$ are relatively prime integers
satisfying the equation $\frac{q_2}{p_2}=\frac{kq_1}{lp_1}$ and $P_1$ is any closed proper Reinhardt subset of $\mathbb C$ (then, obviously,
$P_2$ is of the form $\{\mu_1\zeta^kB(0):\ \zeta\in P_1\}$).

\item Up to a permutation of the components of $f$ and the variables,
the map $f$ has the form
\begin{equation} f(z,w)=((e^{it_1}z^{a_1}+s)^{a_2},
e^{it_2}\exp(2\bar se^{it_1}z^{a_1}+|s|^2)^{-c_2}w^{c_1c_2}),
\end{equation} where $a_1,a_2,c_1,c_2\in\mathbb N,$ $s\in\mathbb C_*$ and $t_1,t_2\in\mathbb R.$ In this case domains have the forms
\begin{align} D_i=\{(z,w)\in\mathbb C^2:\ |w|<C_i\exp(-E_i|z|^{k_i})\},\ i=1,2, \end{align} where $k_1=2a_1,\ k_2=2/a_2$ and $C_1,C_2,E_1,E_2>0.$
\end{enumerate}
\end{theorem}

As mentioned before, in Section~\ref{d=2} we shall also obtain some results related to proper mappings $f:D\to G$ in the case when
$d(D)=d(G)=2.$ It is clear that for any pseudoconvex domain $D$ in $\mathbb C^2,$ $d(D)=2$ if and only if $\log D=\mathbb R^2.$

\section{Proofs}

We start with the following
\begin{lemma}\label{lemwlasciwe}
Let $\varphi:D_1\to D_2$ be a proper holomorphic mapping, where $D_1,D_2\subset\mathbb C^n$ are pseudoconvex Reinhardt domains.
\begin{enumerate}
\item[(a)] Assume that $d(D_2)=0$ and suppose that there is a non-constant holomorphic mapping $\psi:\mathbb C\to D_1.$ Then
$\varphi(\psi(\mathbb C))\subset M.$ \item[(b)] If $\tilde\psi:\mathbb C\to D_2$ is a non constant holomorphic mapping and $d(D_1)=0,$ then
$\varphi^{-1}(\tilde\psi(\mathbb C))\subset M.$
\end{enumerate}
\end{lemma}
\begin{proof}
a) By Lemma 6 in \cite{jar85} there is a nonempty open set $U\subset \mathbb R^n$ and there is a positive $R$ such that for any $v\in U$ the set
$\log D_2$ is contained in $\{x\in\mathbb R^n:\ x_1v_1+\ldots+x_nv_n <R\}.$ Thus, there are linearly independent
$\alpha^1,\ldots,\alpha^n\in\mathbb R^n,$ $\alpha^{\iota}=(\alpha^{\iota}_1,\ldots,\alpha_n^{\iota}),$ such that $D_2$ is contained in
$\{z\in\mathbb C^n:\ |z^{\alpha^{\iota}}|<e^R\},\ \iota=1,\ldots,n.$ Put
\begin{equation}\label{defu}u_{\iota}(z)=|\varphi(\psi(z))^{\alpha^{\iota}}|,\quad
z\in\mathbb C,\quad \iota=1,\ldots,n.
\end{equation}
Obviously $u_{\iota}$ are bounded and subharmonic functions on $\mathbb C,$ so they are constant, say $u_{\iota}=\rho_{\iota},\
\iota=1,\ldots,n.$ It suffices to notice that $\rho_{\iota}=0$ for some $\iota.$ Indeed, if $\rho_{\iota}\neq0$ for every $\iota=1,\ldots,n,$
then obviously $\sum_{j=1}^n\alpha_j^{\iota}\log|\varphi_j(\psi(z))|=\log\rho_{\iota}.$ Applying Cramer rules we would find that the mapping
$\varphi\circ\psi$ would be constant (recall that $\alpha^1,\ldots,\alpha^n$ are linearly independent). However, it would be obviously in
contradiction with the properness of the mapping $\varphi$ (as the mapping $\psi$ is unbounded).

b) Let $\alpha=(\alpha_1,\ldots,\alpha_n)\in\mathbb R_*^n$ and $R>1$ be such that $\log D_1$ is contained in $\{x\in\mathbb R^n:\
x_1\alpha_1+\ldots+x_n\alpha_n<R\}$ and for any $t\in\mathbb R$ the set $\{x\in\mathbb R^n:\ x_1\alpha_1+\ldots+x_n\alpha_n=t\}\cap\log D_1$ is
bounded. Put $u(z)=|z_1|^{\alpha_1}\ldots|z_n|^{\alpha_n}$ for $z\in D_1$. It is well known that the function
\begin{equation}\label{defv}v(z)=v_{\alpha}(z)=\max u(\varphi^{-1}(\tilde\psi(z))),\quad z\in\mathbb C,
\end{equation}
is subharmonic. As it is bounded we find that $v$ is constant. Let $\rho$ be such that $v=\rho.$ Similarly as in the previous part of the proof
it is sufficient to show that $\rho$ is equal to $0.$

Suppose not. One can see that there is a sequence $\{w_{\mu}\}_{\mu=1}^{\infty}\subset D_1$ such that $|w_{\mu}^{\alpha}|=\rho$ for any
$\mu\in\mathbb N$ and $w_{\mu}\to w_0\in\partial D_1,\ (\mu\to\infty).$ Moreover, $|w^{\alpha}|\leq\rho$ for every $w$ such that $\varphi(w)\in
\tilde\psi(\mathbb C).$ Take the supporting hyperplane $H$ of $\log D_1$ at the point $\log w_0$ and let $\beta\in\mathbb R^n$ be such that
$H=\{x\in\mathbb R^n:\ x_1\beta_1+\ldots+x_n\beta_n=\hat\rho\}$ for some $\hat\rho\in\mathbb R.$ Repeating the above reasoning (here the
assumption of the boundedness of $H\cap\log D_1$ is unnecessary) applied to a function $v=v_{\beta}$ (see (\ref{defv})), we find that there is
$\tilde\rho<e^{\hat\rho}$ such that $|w^{\beta}|\leq\tilde\rho$ for any $w\in \varphi^{-1}(\tilde\psi(\mathbb C)).$ However,
$|w_{\mu}^{\beta}|\to e^{\hat\rho},\ (\mu\to\infty),$ which immediately gives a desired contradiction.
\end{proof}

\begin{corollary}\label{theorem3} Let $D,G\subset\mathbb C^n$ be pseudoconvex Reinhardt
domains such that $d(D)=0$ and $d(G)\geq1.$ Then the sets $\prop(D,G)$ and $\prop(G,D)$ are empty.
\end{corollary}
\begin{proof}Take $\alpha=(\alpha_1,\ldots,\alpha_n),\beta=(\beta_1,\ldots,\beta_n)\in\mathbb R^n$ such that $\alpha\mathbb R+\beta\subset\log G.$
Note that for any $z\in G$ the set $\psi_z(\mathbb C)$ is contained in $G,$ where $\psi_z$ is given by a formula
\begin{equation}\psi_z(\zeta)=(z_1e^{\alpha_1\zeta},\ldots,z_ne^{\alpha_n\zeta}),\quad \zeta\in\mathbb C.
\end{equation}

Thus, if $f:D\to G$ (or $g:G\to D$) would be a proper holomorphic mapping then, by Lemma \ref{lemwlasciwe} $f^{-1}(G)\subset M$ (resp.
$g(G)\subset M$). This immediately gives a contradiction.
\end{proof}

\begin{lemma}\label{brod}Let $D,G\subset\mathbb C^n$ be domains. Assume that $D$ is bounded and $G$ is not Brody-hyperbolic. Then there is
no proper holomorphic mapping from $D$ onto $G.$
\end{lemma}

\begin{proof} Suppose that $\varphi:D\to G$ is a proper holomorphic
mapping. Put $A=\{z\in D: \det\varphi'(z)=0\}.$ The set $A$ is a variety in $D$ and, by the properness of $\varphi,$ $A\neq D.$ Moreover, there
is an integer $m$ such that $\#\varphi^{-1}(w)=m$ for any $w\in G\setminus\varphi(A).$

Put
$$\pi_k(\lambda)=\sum_{1\leq
i_1<\ldots<i_k\leq m}\lambda_{i_1}\ldots\lambda_{i_k},\quad \lambda=(\lambda_1,\ldots,\lambda_m)\in\mathbb C^m,\ k=1,\ldots,m$$ and
$\pi=(\pi_1,\ldots,\pi_m).$ Moreover, for $z^j=(z^j_1,\ldots,z^j_n)\in\mathbb C^n,\ j=1,\ldots,m,$ define
\begin{equation}\label{sigma}
\sigma(z^1,\ldots,z^m):=(\pi(z_1^1,\ldots,z_1^m),\ldots,\pi(z_n^1,\ldots,z_n^m)).
\end{equation} Obviously $\sigma:(\mathbb C^n)^m\to\mathbb C^{nm}$
is a proper holomorphic mapping with multiplicity equal to $(m!)^n.$

Let $\varphi^{-1}(w)=\{\zeta_1(w),\ldots,\zeta_m(w)\},\ w\in G\setminus\varphi(A).$ Since $\varphi$ is locally biholomorphic near any
$\zeta_i(w),$ $i=1,\ldots,m,$ and the mapping $\sigma$ given by the formula (\ref{sigma}) is symmetric, we find that the mapping
$\psi=\sigma\circ(\zeta_1,\ldots,\zeta_m)$ is holomorphic in $G\setminus\varphi(A).$ Because $\varphi(A)$ is analytic in $G$ and $\psi$ is
bounded, we may extend $\psi$ to bounded holomorphic mapping on the whole $G.$ Let $\tilde\psi$ be such an extension. Take any $\gamma:\mathbb
C\to G$ non-constant and holomorphic. Then $\tilde\psi\circ\gamma$ is bounded and holomorphic on $\mathbb C;$ in particular,
$\tilde\psi\circ\gamma$ is constant.

Let us take any $z'\in\mathbb C.$ If $\gamma(z')$ belongs to $G\setminus\varphi(A),$ then obviously
$\tilde\psi(\gamma(z'))=\sigma(\zeta_1(\gamma(z)),\ldots,\zeta_m(\gamma(z))).$ Suppose now that $\gamma(z')$ is a critical value of $\varphi.$
Let $x=(x_1,\ldots,x_m)\in(\mathbb C^n)^m$ be such that $\tilde\psi(\gamma(z))=\sigma(x),\ z\in\mathbb C.$

Take any $\zeta$ such that $\varphi(\zeta)=\gamma(z'),$ and let $(\zeta_n)\subset D\setminus A$ be such that $\zeta_n\to \zeta.$ Observe that
$\sigma(\varphi^{-1}(\varphi(\zeta_n)))=\tilde{\psi}(\varphi(\zeta_n))\to\sigma(x).$ In particular, using properness of $\sigma,$ we find that
$\zeta\in\sigma^{-1}(\sigma(x)),$ so we have shown that $\varphi^{-1}(\gamma(z'))\subset\sigma^{-1}(\sigma(x_1,\ldots,x_1)).$

It follows that for any $w\in\gamma(\mathbb C)$ $\varphi^{-1}(w)$ is contained in the finite set $\sigma^{-1}(\sigma(x)).$ Since the mapping
$\gamma$ is unbounded, we immediately get a contradiction with the properness of $\varphi.$
\end{proof}

\begin{remark}
Since the mapping $\mathbb C\setminus\{0,1\}\ni z\to \frac{1}{z(z-1)}\in\mathbb C$ is proper, the above theorem is not true if we only assume
that the domain $D$ is Brody-hyperbolic (instead of bounded). On the other hand, since in the class of pseudoconvex Reinhardt domains the
property of Brody-hyperbolicity means, up to algebraic mappings, the boundedness, we easily see that there is no proper holomorphic mapping
between hyperbolic and non-hyperbolic pseudoconvex Reinhardt domains.
\end{remark}

For a Reinhardt domain $D$ in $\mathbb C^n$ let $I(D)$ denote the set of $i=1,\ldots,n$ for which the intersection $ V_i\cap D$ is not
$c-$hyperbolic (viewed as a domain in $\mathbb C^{n-1}$). Put
\begin{equation}D^{hyp}=D\setminus (\bigcup_{i\in I(D)} V_i).
\end{equation} It is
clear that $D^{hyp}=D$ if $D$ is $c$-hyperbolic or $D\subset\mathbb C_*^n.$ In the sequel by $\widehat D^{hyp}$ we shall denote the set
$(\widehat D)^{hyp}.$

Now we are in position to formulate the following

\begin{theorem}\label{theorem2}
Let $D_1,D_2$ be pseudoconvex Reinhardt domains in $\mathbb C^2$ such that $\log D_i$ contains no affine lines, $i=1,2$ (i.e.
$d(D_1)=d(D_2)=0$). If $\varphi:D_1\to D_2$ is a proper holomorphic mapping, then $\varphi(D_1^{hyp})\subset D_2^{hyp}$ and the restriction
$\varphi|_{D_1^{hyp}}:D_1^{hyp}\to D_2^{hyp}$ is proper.
\end{theorem}
\begin{proof}Obviously it suffices to prove the following
statement:
\begin{enumerate}
\item[1)] If $D_1\cap V_1\in\{\mathbb C,\mathbb C_*\}$ then $\varphi(D_1\cap V_1)$ is contained either in $V_1$ or in $V_2.$
In particular, if moreover $D_2\cap V_2$ is c-hyperbolic or empty, then $\varphi(D_1\cap V_1)\subset V_1.$ \item[2)] If $D_2\cap V_1\in\{\mathbb
C,\mathbb C_*\}$ and $D_1\cap V_2$ is neither $\mathbb C$ nor $\mathbb C_*,$ then $D_1\cap V_1\in\{\mathbb C,\mathbb C_*\}$ and
$\varphi^{-1}(D_2\cap V_1)\subset V_1.$
\end{enumerate}
1) From Lemma \ref{lemwlasciwe}(a) applied to the mapping $\psi(z)=(0,e^z),\ z\in\mathbb C,$ we get that $\varphi(D_1\cap V_1)\subset M.$ It
follows that $\varphi_1(0,z)\varphi_2(0,z)=0$ for any $z\in\mathbb C_*,$ so $\varphi_1(0,\cdot)\equiv0$ or $\varphi_2(0,\cdot)\equiv0.$ The
second statement is clear.

2) If $D_1\cap V_1$ were neither $\mathbb C$ nor $\mathbb C_*,$ then $D_1$ would be biholomorphic to a bounded domain, which obviously
contradicts Lemma \ref{brod}. Thus $D_1\cap V_1\in\{\mathbb C,\mathbb C_*\}.$

Suppose that $D_1\cap V_2$ is non-empty (in the other case Lemma~\ref{lemwlasciwe}~(b) finishes the proof). Pseudoconvexity implies that
$\pi_1(D_1)$ is a bounded subset of $\mathbb C,$ where $\pi_1:\mathbb C^2\to\mathbb C$ denotes a projection onto the first variable. Thus, a
function given by the formula
\begin{equation}v(z):=\max|\pi_1(\varphi^{-1}(0,z))|,\quad
z\in\mathbb C_*
\end{equation} is constant (as bounded and subharmonic).
Moreover, from Lemma \ref{lemwlasciwe} we get that $\varphi^{-1}(D_2\cap V_1)\subset M.$

Now, one may easily verify that $v=0.$
\end{proof}

\begin{corollary}\label{niezmiennik}
Let $D_1,D_2\subset\mathbb C^2$ be Reinhardt domains such that $d(D_1)=d(D_2)=0.$ Assume that $\prop(D_1,D_2)$ is non-empty. Then
\begin{equation}(s(D_1),s_*(D_1),t(D_1))=(s(D_2),s_*(D_2),t(D_2)).\end{equation}
\end{corollary}
\begin{proof}Any proper holomorphic map between domains $D_1,D_2$ may be extended to the proper map between their envelopes of holomorphy
$\widehat D_1,\widehat D_2,$ respectively. Moreover, it is well known (see Corollary~0.3~\cite{isaev}) that if there exists a proper holomorphic
mapping between two given bounded domains, then there also exists an elementary algebraic proper holomorphic mapping between these domains.

Thus, our result is a direct consequence of Theorem~\ref{theorem2} and properties of algebraic mappings.
\end{proof}

Note, that applying the methods used in previous theorems we may easily show

\begin{proposition}\label{aa} There are no proper holomorphic mappings between domains $D,G$ and $G,D$ in the following cases:
\begin{enumerate}[1.]
\item $D\subset\{z\in D':\ u(z)<0\},$ where $u$ is some plurisubharmonic function on the domain $D'\subset\mathbb C^n,$ non constant on $D$ and $G=\mathbb C^n\setminus E$ for some pluripolar set $E$ in $\mathbb C^n.$
\item $D$ is hyperconvex (i.e. there is a negative plurisubharmonic exhaustion function for $D$) and $G$ is not
Brody-hyperbolic.
\end{enumerate}
\end{proposition}
\begin{proof}1. Obviously $\prop(\mathbb C^n\setminus E,D)=\emptyset.$

Suppose that $\varphi:D\to \mathbb C^n\setminus E$ is proper and holomorphic. Put $v(z)=\max u(\varphi^{-1}(z)),$ $z\in\mathbb C^n\setminus E.$
It is seen that the function $v$ is constant. In particular, there is $\rho<0$ such that $u\leq\rho$ and $u(w_0)=\rho$ for some $w_0\in D;$ a
contradiction.

2. It is clear that the set $\prop(G,D)$ is empty. Suppose that $\varphi:D\to G$ is proper and holomorphic. Let $u$ be a negative
plurisubharmonic exhaustion function for $D$ and let $\psi:\mathbb C\to G$ be a non-constant holomorphic mapping. Put $v(\zeta)=\max
u(\varphi^{-1}(\psi(\zeta))),\ \zeta\in\mathbb C.$ The function $v$ is subharmonic on $\mathbb C.$ Since $v<0,$ it is constant. So we find that
$\varphi^{-1}(\psi(\mathbb C))$ is a relatively-compact subset of $D;$ a contradiction.
\end{proof}

\begin{proposition}\label{dim}
Let $D,G\subset\mathbb C^2$ be pseudoconvex Reinhardt domains. If $d(D)\neq d(G),$ then there is no proper holomorphic mapping between $D$ and
$G.$
\end{proposition}
\begin{proof}
In the view of Corollary~\ref{theorem3} it suffices to consider the case when $d(D)=2$ or $d(G)=2.$ However, this immediately follows from
Proposition~\ref{aa}.
\end{proof}

\begin{proof}[Proof of Theorem~\ref{niezmienniki}] It is a direct consequence of Corollary~\ref{niezmiennik} and Proposition~\ref{dim}.
\end{proof}

\begin{proof}[Proof of Theorem~\ref{main}]
Take any $f:D_1\to D_2$ proper and holomorphic and suppose that it is non-elementary. Let $f:\widehat D_1\to \widehat D_2$ also denotes its
extension to the proper mapping between the envelopes of holomorphy of $D_1$ and $D_2.$ By Proposition~\ref{dim} $d(D_1)=d(D_2).$

First we consider the case $d(D_1)=d(D_2)=0.$ Then, by Theorem~\ref{theorem2} the restriction $f|_{\widehat D_1^{hyp}}:\widehat D_1^{hyp}\to
\widehat D_2^{hyp}$ is also proper.

If $s(D_1)=2$ or $s_*(D_1)=t(D_1)=1,$ then $\widehat D_1^{hyp}$ would be contained in $\mathbb C_*^2$ and from the description obtained in
\cite{isaev} we find that $f|_{\widehat D_1^{hyp}}$ would be elementary algebraic. It is clear that the identity principle gives a
contradiction.

Therefore, we may assume, that $t(D_1)=t(D_2)=2,\ s(D_1)=s(D_2)=1$ and $s_*(D_1)=s_*(D_2)=0.$ Up to a permutation of components we may suppose
that $\widehat D_i\cap V_2=V_2$ and $\widehat D_i\cap V_1$ is bounded, $i=1,2.$ Therefore, there are $k_1,k_2\in\mathbb N$ such that $\widehat2
D_i$ is contained in $\{(z,w)\in\mathbb C^2:\ |z||w|^k<c_i\}$ for some positive constants $c_i,\ i=1,2.$ It follows that $\Phi_{A_i},$ where
$A_i=\left(
\begin{array}{cc}
  1 & k_i \\
  0 & 1 \\
\end{array}
\right),$ is a biholomorphic mapping from $\widehat D_i^{hyp}$ onto the bounded set $\Phi_{A_i}(\widehat D_i^{hyp}),\ i=1,2.$ In particular,
\begin{equation}g:=\Phi_{A_2}\circ f\circ \Phi_{A_1^{-1}}:\Phi_{A_1}(\widehat D_1^{hyp})\to\Phi_{A_2}(\widehat D_2^{hyp})\end{equation} is a proper
holomorphic mapping between two bounded domains in $\mathbb C^2.$ Now, using description obtained in \cite{isaev} it is straightforward to
observe that two possibilities may hold:
\begin{enumerate}\item[(i)]\label{ok} $\widehat
D_i^{hyp} =\{(z,w)\in\mathbb C^2:\ C_i|z|^{p_i}|w|^{p_ik_i+q_i}<1,\ 0<|w|<C'_i\},$ where $p_i,q_i$ are relatively prime integers such that
$p_ik_i+q_i>0,\ p_i>0,\ q_i\leq0$ and $C_i,C_i'>0,\ i=1,2.$ \item[(ii)]\label{nok} $\widehat D_1^{hyp}=\{(z,w)\in\mathbb C^2:\
0<|w|<C_1\exp(-E_1|z|^{2a_1}|w|^{2k_1a_1-2b_1}),$ and $\widehat D_2^c=\{(z,w)\in\mathbb C^2:\
0<|w|<C_2\exp(-E_2|z|^{2/a_2}|w|^{2k_2/a_2-2b_2/a_2c_2}),$ where $a_i,b_i,c_i\in\mathbb N,\ C_i,E_i>0,\ i=1,2.$
\end{enumerate}

First suppose that (i) holds. From \cite{isaev} it follows that $g$ must be of the form $g(z,w)=(\lambda_1 z^aw^b B(C_1z^{p_1}w^{q_1}),\lambda_2
w^c),\ (z,w)\in\Phi(\widehat D_1^{hyp}),$ where $a,b,c\in\mathbb Z,\ a,c>0, aq_1-bp_1<0,\ \frac{q_2}{p_2}=\frac{aq_1-bp_1}{cp_1},$ $B$ is a
Blaschke product non-vanishing at $0$ and $\lambda_1,\lambda_2\in\mathbb C_*.$ Put $\tilde q_i=p_ik_i+q_i.$ It is obvious that $p_i$ and $\tilde
q_i$ are relatively prime. Moreover, from the form of $\widehat D_i$ we get that $\tilde q_i>0.$ An easy computation gives
$$f(z,w)=(\mu_1z^aw^{ak_1-ck_2+b}B(C_1z^{p_1}w^{\tilde
q_1}),\mu_2w^c),\quad (z,w)\in\widehat D_1^{hyp},$$ for some constants $\mu_1,\mu_2.$ Since $f$ may be extended properly on $\widehat D_1,\
ak_1-ck_2+b=0.$ Moreover, it is clear that $\frac{\tilde q_2}{p_2}=\frac{a\tilde q_1}{cp_1}.$

It is straightforward to see that any Reinhardt subdomain of $\widehat D_1$ mapped properly by $f$ onto a Reinhardt domain and whose envelopes
of holomorphy coincides with $\widehat D_1$ is equal to $\widehat D_1\setminus P_1\times \{0\},$ where $P_1$ is any closed Reinhardt subset of
$\mathbb C.$

Now suppose that (ii) holds. Denote $m_1:=k_1a_1-b_1,\ m_2:=k_2c_2-b_2.$ Similarly as before, taking into account the form of $\widehat D_1$ and
$\widehat D_2$ one can see that $m_1,m_2\geq0.$

For $s\in\mathbb C_*$ and $t_1,t_2\in\mathbb R$ put $h_1(z):=e^{it_1}z+s,\ h_2(z)=e^{it_2}\exp(2\bar se^{it_1 z}+|s|^2),\ z\in\mathbb C.$ An
easy calculation and formula for the mapping $g$ (see \cite{isaev}) give
\begin{align}f(z,w)=(h_1(z^{a_1}w^{m_1})^{a_2}h_2(z^{a_1}w^{m_1})^{m_2}
w^{-m_2c_2},h_2(z^{a_1}w^{m_2})^{-c_2}w^{c_1c_2}),\\\nonumber (z,w)\in\widehat D_1^{hyp}.\end{align} Since $f$ may be extended through $V_1,$
$m_1=m_2=0.$

Finally, one may easily verify that any Reinhardt subdomain of $\widehat D_1$ whose envelope of holomorphy coincides with $\widehat D_1$ and
which is mapped properly by $f$ onto a Reinhardt domain is equal to $\widehat D_1.$
\end{proof}

\section{Remarks on the proper holomorphic mappings $f:D\to G$ between Reinhardt domains in case $d(D)=d(G)=2$}\label{d=2}

It is well known that the structure of $\aut(\mathbb C^2),\ \aut(\mathbb C_*^2),$ $\aut(\mathbb C\times\mathbb C_*)$ is very complicated and the
full description of these groups seems to be not known. Proper maps are harder to deal with, so description of the set of proper holomorphic
mappings between pseudoconvex Reinhardt domains $D_1$ and $D_2$ in the case, when $\log D_i=\mathbb R^2,\ i=1,2,$ is more difficult.

Below we present some partial results related to these problems.

\begin{proposition}\label{propozycja}
The sets $\prop(\mathbb C\times\mathbb C,\mathbb C\times\mathbb C_*),\ \prop(\mathbb C\times\mathbb C,\mathbb C_*\times\mathbb C_*)$ and
$\prop(\mathbb C\times\mathbb C_*,\mathbb C_*\times\mathbb C_*)$ are empty.
\end{proposition}

\begin{proof} First suppose that $f:\mathbb C^2\to\mathbb C\times\mathbb C_*$ is proper and holomorphic. Obviously there exists a holomorphic mapping
$\psi:\mathbb C^2\to\mathbb C^2$ such that $f=(\psi_1,e^{\psi_2}).$ One can easily verify that the mapping $\psi$ is proper; in particular
$\psi$ is surjective. Thus there is a discrete sequence $(z_n)_{n\in\mathbb N}\subset\mathbb C^2$ such that $\psi(z_n)=(0,2n\pi i),\ n\in\mathbb
N.$ It follows that $f(z_n)=(0,1)$ for $n\in\mathbb N.$ From this we immediately get a contradiction.

To show that $\prop(\mathbb C^2,\mathbb C_*^2)=\emptyset$ we proceed similarly.

Now suppose that $g:\mathbb C\times\mathbb C_*\to\mathbb C_*^2$ is holomorphic and proper. It is seen that there exists a holomorphic mapping
$\varphi:\mathbb C^2\to\mathbb C^2$ such that $g(z,e^w)=(e^{\varphi_1(z,w)},e^{\varphi_2(z,w)})$ for $z,w\in\mathbb C.$

Fix $z\in\mathbb C$ and put $\tilde g_i=g_i(z,\cdot),\ \tilde\varphi_i=\varphi_i(z,\cdot),\ i=1,2.$ Since $\tilde
g_i(e^w)=e^{\tilde\varphi_i(w)},$ we find that $\tilde\varphi_i'(w)=\zeta_i(e^w),\ w\in\mathbb C,$ where $\zeta_i$ is a holomorphic function
given by the formula $\zeta_i(\lambda)=\frac{\lambda\tilde g_i'(\lambda)}{\tilde g_i(\lambda)},\ \lambda\in\mathbb C_*.$ Expanding $\zeta_i$ to
the Laurent series gives $\tilde\varphi_i(w)=a_iw+\sum_{n\in\mathbb Z_*}a_{in}e^{nw}$ for some $a_i=a_i(z),a_{in}=a_{in}(z)\in\mathbb C.$

Thus, there is a holomorphic mapping $\hat\varphi_i(\cdot)=\hat\varphi_i(z,\cdot)$ on $\mathbb C_*$ such that
$\tilde\varphi_i(w)=a_iw+\hat\varphi_i(e^w),\ w\in\mathbb C.$ Since $e^{a_iw}=\frac{\tilde g_i(e^w)}{e^{\hat\varphi_i(e^w)}}$ we immediately
find that $a_i\in\mathbb Z,\ i=1,2.$

Therefore $\varphi_i(z,w)=a_i(z)w+\hat\varphi_i(z,e^w),\ z,w\in\mathbb C,\ i=1,2.$ In particular
\begin{equation}g(z,w)=(w^{a_1(z)}e^{\hat\varphi_1(z,w)},w^{a_2(z)}e^{\hat\varphi_2(z,w)}),\quad
(z,w)\in\mathbb C\times\mathbb C_*.\end{equation}

It is straightforward to verify that $a_i(z)=\frac{1}{2\pi i}\int_{\partial \mathbb D}\frac{\frac{\partial g_i}{\partial
\lambda}(z,\lambda)}{g_i(z,\lambda)}d\lambda,\ z\in\mathbb C,$ whence $a_i$ is constant (recall that $a_i(z)\in\mathbb Z$) and therefore
$\hat\varphi_i$ is holomorphic on $\mathbb C\times\mathbb C_*,\ i=1,2$.

Note that we may assume that $a_2=0$ (if $a_1a_2\neq0$ one may compose $g$ with a proper holomorphic mapping $F:\mathbb C_*^2\to\mathbb C_*^2$
given by the formula $F(z,w)=(z^{a_2},\frac{w^{a_1}}{z^{a_2}})$).

Put $$h(z,w)=(w^{a_1}e^{\hat\varphi_1(z,w)},\hat\varphi_2(z,w)),\quad (z,w)\in\mathbb C\times \mathbb C_*,$$ and notice that the mapping
$h:\mathbb C\times\mathbb C_*\to\mathbb C_*\times\mathbb C$ is proper.

Now, in order to get a contradiction one may proceed exactly as in the case of $\prop(\mathbb C^2,\mathbb C\times\mathbb C_*).$
\end{proof}

\begin{corollary}
$\prop(A\times\mathbb C,A\times\mathbb C_*)$ is empty for any domain $A\subset\mathbb C.$
\end{corollary}

\begin{proof}If $\#(\mathbb C\setminus A)\leq 1$ the result follows directly from Proposition~\ref{propozycja}. Assume that
$\#(\mathbb C\setminus A)>1$ and let $f:A\times\mathbb C\to A\times\mathbb C_*$ be proper and holomorphic. By the uniformization theorem there
is an universal covering $\pi:\mathbb D\to A$ and there is $\psi\in\mathcal O(\mathbb D\times\mathbb C,\mathbb D)$ such that
$$f(\pi(\lambda),w)=(\pi(\psi(\lambda,w)),f_2(\pi(\lambda),w))\quad \text{for any}\quad (\lambda,w)\in\mathbb D\times\mathbb C.$$
Fix any $\lambda\in\mathbb D$ and note that the mapping $\psi(\lambda,\cdot)$ is constant. From properness of $f$ it easily follows that the
mapping $f_2(\pi(\lambda),\cdot):\mathbb C\to\mathbb C_*$ is proper; a contradiction.
\end{proof}

\begin{remark}Since $\phi:\mathbb C_*\ni z\to z+1/z\in\mathbb C$ is proper, there exist proper holomorphic maps from $\mathbb C_*^2$
onto $\mathbb C^2,$ from $\mathbb C_*^2$ onto $\mathbb C\times\mathbb C_*,$ and from $\mathbb C\times\mathbb C_*$ onto $\mathbb C^2.$ Obviously
such maps cannot be algebraic.

On the other hand, the above results and the ones obtained in \cite{isaev} and \cite{ja} imply that if there exists a proper holomorphic mapping
between two Reinhardt domains $D_1,D_2\subset\mathbb C^2$ such that $\alpha\mathbb R+\beta$ is not contained in $\log D_1$  for any
$\alpha\in\mathbb Q^2,\ \beta\in\mathbb R^2$ (hence also in $\log D_2$, see \cite{ja}), then there also exists an elementary algebraic mapping
between these domains.
\end{remark}

\end{document}